# Trisections and Totally Real Origami

Roger C. Alperin

## 1 CONSTRUCTIONS IN GEOMETRY.

The study of methods that accomplish trisections is vast and extends back in time approximately 2300 years. My own favorite method of trisection from the Ancients is due to Archimedes, who performed a "neusis" between a circle and line. Basically a nuesis (or use of a marked ruler) allows the marking of points on constructed objects of unit distance apart using a ruler placed so that it passes through some known (constructed) point $P$.

Here is Archimedes' trisection method: Given an acute angle between rays $r$ and $s$ meeting at the point $O$, construct a circle $K$ of radius one at $O$, and then extend $r$ to produce a line that includes a diameter of $K$. The circle $K$ meets the ray $s$ at a point $P$. Now place a ruler through $P$ with the unit distance $CD$ lying on the circle at $C$ and diameter at $D$, on the opposite ray to $r$. That the angle $\angle ODP$ is the desired trisection is easy to check using the isosceles triangles $\triangle DCO$ and $\triangle COP$ and the exterior angle of the triangle $\triangle PDO$. As one sees when trying this for oneself, there is a bit of "fiddling" required to make everything line up as desired; that fiddling is also essential when one does origami.

The Greeks' use of neusis gave us methods not only for trisections of angles but also for extraction of cube roots. Thus in a sense all cubics could be solved by the Greeks using geometric methods. These methods were highly developed by Archimedes and Apollonius and then further by Diocles and Nicomedes [10]. The latter found a trisection method using a neusis between lines.

Figure 1: Archimedes trisection

A problem, reported by Pappus, that comes from one of Apollonius' now lost works is to mark a fixed distance, say a unit length, as a chord of a given conic so that the extension of the chord passes through a given point $P$. Apollonius tried to reduce this and other neusis problems to constructions with ruler and compass; sometimes he was successful. However, the chord construction just described in general requires the solution of a polynomial equation $g(x^2) = 0$, where $g$ is a quartic; this ensures that the equation has a solvable Galois group, but generally its solutions are not constructible by ruler and compass. Other neusis constructions between conics do not usually have solvable Galois group; indeed, they typically entail solution of irreducible equations having degree sixteen. Some marked ruler constructions that are nonsolvable have been discussed recently in [4].

Neusis constructions between lines are called Viéte marked ruler constructions [13]. The set of all real numbers that can be obtained via Viéte marked ruler constructions is a subfield $\mathbb{V}$ of $\mathbb{R}$. Its elements are called Viétens. The set of complex numbers that can be obtained by such constructions is $\mathbb{V} \oplus i\mathbb{V}$, the subfield of $\mathbb{C}$ consisting of all algebraic numbers that lie in normal extensions of the rationals with degrees of the form $2^a 3^b$. The field $\mathbb{V} \oplus i\mathbb{V}$ has also been realized as the field of Mira constructions [5], the field of intersections of constructible conics [12], and the field of origami numbers [1].

Our aim here is to characterize the totally real Viétens, thus providing the analogue of Hilbert's characterization of totally real Euclidean numbers (i.e., the Pythagorean numbers). Axioms and a characterization of these numbers are given in terms of the origami numbers. The treatment here is mostly self-contained, although we appeal to important facts from field theory. In extending from the Pythagorean numbers to the totally real Viétens, it is important to understand the constructions involving trisections. Gleason's article [6] has been an inspiration in this study. We introduce the origami axiom (T), which allows certain restricted foldings. It was motivated by Abe's method described in [8]. (In fact, Abe's method has a classical counterpart using the trisectrix of Maclaurin.) We show that axiom (T) is just the right axiom for allowing us to construct precisely the totally real Viétens. In the last section, we discuss the field of constructions studied by Gleason and relate the solution of "Alhazen's Problem" to the corresponding field theory.

## 2 TOTALLY REAL ORIGAMI.

The field $\mathbb{E}$ of *Euclidean numbers* comprises precisely those points in the complex numbers $\mathbb{C}$ that can be constructed by means of the compass and straightedge, starting with a marked origin and unit length. This field can also be described as the set of algebraic numbers that lie in normal (Galois) extensions of the rationals $\mathbb{Q}$ whose degrees are powers of two. Recall that an algebraic number $\alpha$ is a complex number that satisfies an irreducible polynomial with rational coeffficients. The roots of this polynomial are called the (algebraic) conjugates of $\alpha$.

In his simplification of the foundations of geometry [7], Hilbert described the field $\mathbb{P}$ of *Pythagorean numbers* which is the totally real subfield of the real Euclidean numbers (i.e, all algebraic conjugates are real). The field $\mathbb{P} \oplus i\mathbb{P}$ is the set of complex numbers that can be constructed using only the straightedge and dividers (called a *scale* by Hilbert). A pair of dividers, which is a weaker tool than a compass, enables one to transfer distances. Constructions with dividers give only those lines whose slopes are algebraic real numbers having only real conjugates. For example, the real number $\sqrt{1+\sqrt{2}}$ is not totally real (it has conjugate $\sqrt{1-\sqrt{2}}$ ) and consequently not all right triangles constructible in Euclidean geometry can be constructed in Pythagorean geometry. However, any regular polygon that is Euclidean constructible is also Pythagorean constructible. We now develop some of the key ideas relating totally real extensions and Pythagorean constructions.

Of central importance in understanding the Pythagorean constructions is the Artin-Schreier view of totally real numbers [9 (vol. 2), Theorem 11.5]. The Artin-Schreier theorem implies that in a subfield $F$ of the real algebraic numbers an element $\alpha$ is totally positive (i.e., all algebraic conjugates of $\alpha$ are positive) if and only if $\alpha$ is a sum of squares in $F$.

For example, a quadratic extension of a totally real field $K$ is obtained by adjoining to $K$ a root $\sqrt{a}$, $a \in K$ of a quadratic (reduced) irreducible polynomial $x^2 - a$. An element $\alpha$ of $K(\sqrt{a})$ has the form $\alpha = u + v\sqrt{a}$ and an algebraic conjugate of $\alpha$ is obtained by applying to $\alpha$ an extension of an injective homomorphism $\sigma : K \to \mathbb{R}$ defined as follows: an extension of $\sigma$ for $\epsilon \in \{\pm 1\}$, $\sigma_\epsilon : K(\sqrt{a}) \to \mathbb{R}$, is defined as $\sigma_\epsilon(u + v\sqrt{a}) = \sigma(u) + \epsilon\sigma(v)\sqrt{a}$, for $u, v \in K$; the map $\sigma_\epsilon$ is also an injective homomorphism. We thus have the following easy result:

**Lemma 2.1.** *The quadratic extension $K(\sqrt{a})$ of the totally real field $K$ is also totally real if and only if the element $a$ of $K$ is totally positive.*

To see why Lemma 2.1 is true observe first that if $a$ in $K$ is totally positive, then it follows from the Artin-Schreier theorem that $a$ is a sum of squares in $K$, say $a = x_1^2 + \cdots + x_n^2$ with $x_1, x_2, \ldots, x_n$ in $K$. Hence the algebraic conjugates of $\sqrt{a}$ are $\pm\sqrt{\sigma(x_1)^2 + \cdots + \sigma(x_n)^2}$, which are real for any embedding $\sigma : K \to \mathbb{R}$. Conversely, if $x = \sqrt{a}$ ($\notin K$) is a totally real number, then $x$ satisfies $x^2 - a$. Arguing by contradiction, we see that if some conjugate $\sigma(a)$ is negative, then $\sigma(x)^2 - \sigma(a)$ is nonzero for an extension of $\sigma$ to $K(\sqrt{a})$. This contradiction shows that all conjugates $\sigma(a)$ are positive, whence $a$ is totally positive. We shall later extend this simple idea to cubics.

The Pythagorean field is built up from the rationals by totally real quadratic extensions. The Artin-Schreier theorem and Lemma 2.1 imply that the totally real subfield of the Euclidean field $\mathbb{E}$ of numbers constructible with ruler and compass is the Pythagorean field. Each quadratic extension obtained by extraction of a square root of a totally positive element $a$ can also be constructed using field operations and angle bisections. The element $\sqrt{a} = \sqrt{x_1^2 + \cdots + x_n^2}$ can be constructed inductively from $x_n$ and $\sqrt{x_1^2 + \cdots + x_{n-1}^2}$ using angle bisections. The initial case $\sqrt{x_1^2 + x_2^2}$ is handled by constructing $x_1$ and $x_2$ on perpendicular axes, thereby producing a right triangle with hypotenuse of the desired length. By bisecting one of the acute angles of the triangle and reflecting the hypotenuse across the bisector, we can mark the desired length on the real axis.

**Origami axioms for Pythagorean numbers.** The basic origami axioms specify an origin and a unit length and the following folding (or construction) rules:

(L) A line connecting two constructible points can be folded.

(P) A point of coincidence of two folded lines is a constructible point.

The additional origami axiom that characterizes the Pythagorean constructions is simply then

(B) The line bisecting any given constructed angle can be folded.

The reader should be aware that we allow as a fold line using a special case of axiom (B), the perpendicular to a given line at a given point. In this way we can also fold the perpendicular bisector of two constructed points $A$ and $B$, modifying slightly the standard Euclidean method (Figure 2). We proceed by constructing (i.e., folding) the perpendiculars at $A$ and $B$ to the line $AB$ and bisect the right angles at $A$ and $B$. The lines of bisection intersect at $C$ and $D$. The line $CD$ is the perpendicular bisector of the segment $AB$. Thus we can also reflect a given point not on a given line across that line. We leave as an exercise the problem of constructing the perpendicular to a given line from a point not on that line. (Hint: Use the midpoint triangle of a triangle with vertex at the given point and opposite side on the given line.)

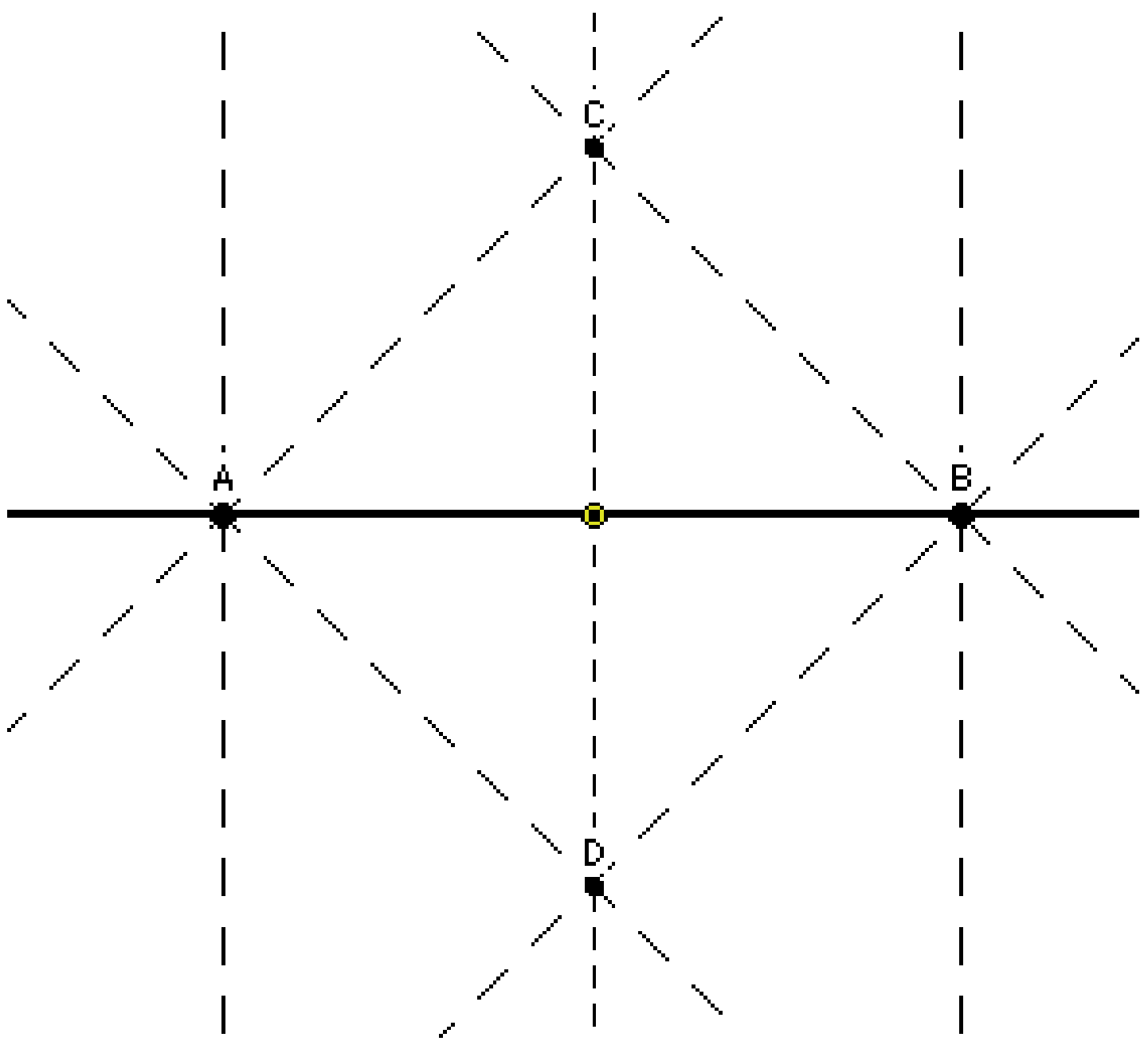

Figure 2: Perpendicular bisection

Once we can construct midpoint bisectors it is easy to see by repeated application of (P) and (L), that the points constructible using the origami axioms form a subfield of the complex numbers and that the slopes of constructible lines form its real subfield. The details can be found in [1]. Basically, the polar representation of complex numbers reduces the field operations to operations on lengths and angles. Many of the proofs use elementary properties of similar triangles. For example, the construction of the reciprocal of a constructed real length $m$ involves the construction a line of slope $m$ and a perpendicular to it, which has slope $-\frac{1}{m}$.

Let $\mathbb{P}$ be the smallest field that contains the rationals and is closed under the operation $x \mapsto \sqrt{1+x^2}$ where $\sqrt{\phantom{x}}$ signifies the positive square root (as noted earlier, $\mathbb{P}$ is the field of Pythagorean numbers). The rules (P), (L) and (B) allow us to construct precisely all the points in $\mathbb{P} \oplus i\mathbb{P}$ where $\mathbb{P}$ is the field of Pythagorean numbers (for more details, see [1]). This simplifies the totally real axioms in [3].

**Theorem 2.2.** *The field $\mathbb{P}$ is the totally real subfield of the field $\mathbb{E}$ of Euclidean numbers. The field $\mathbb{P} \oplus i\mathbb{P}$ is the field of complex numbers constructible by using the axioms (P), (L) and (B).*

**Folding the regular pentagon.** To construct a regular pentagon using the Pythagorean (origami) construction axioms, it is necessary to construct the totally real lengths $s = \sin(2\pi/5)$ and $c = \cos(2\pi/5)$. Starting from the equation

$$\frac{x^5 - 1}{x - 1} = x^4 + x^3 + x^2 + x + 1 = 0$$

for a primitive fifth root of unity $x$, it is easy to see that $t = x + x^{-1} = 2c$ satisfies $t^2 + t - 1 = 0$, so $\cos(2\pi/5) = (\sqrt{5} - 1)/4$. The construction of $\sqrt{5} = \sqrt{1^2 + 2^2}$ can be accomplished with bisections, which implies that $c$ can be constructed by applying the axioms (P), (L), and (B). Since $\cos(\pi/5) = (\sqrt{5} + 1)/4$ we can construct $s$ by Pythagorean constructions using the identity $\sin(2\pi/5)^2 = \cos(\pi/5)^2 + (1/2)^2$.

# 3 ORIGAMI TRISECTOR AXIOM.

We add to the set of axioms (P), (L), and (B) the following "trisection axiom" that allows new fold lines to be created:

(T) Given constructed points $P$ and $Q$ and a constructed line $l$ through $P$ then we can construct (fold) any line $l'$ with the property that the reflection of $Q$ in $l'$ lies on $l$ and the reflection of $P$ in $l'$ lies on the perpendicular bisector of the segment $PQ$.

This axiom is based on the trisection method of Abe discussed in [8]. Our main theorem gives a characterization of the totally real Viétens. For example, $\sqrt[3]{2}$ is an element of $\mathbb{V}$ that does not belong to this totally real field.

The statements of our main results follow. We give their proofs in the following sections.

**Theorem 3.1.** *The field of numbers that can be constructed given an origin and unit length by application of axioms (P), (L), (B), and (T) is $\mathbb{BT} \oplus i\mathbb{BT}$, where $\mathbb{BT}$ is the totally real subfield of $\mathbb{V}$. The field $\mathbb{BT}$ is the smallest field containing the (finite) slopes of all lines that can be folded using bisections and trisections.*

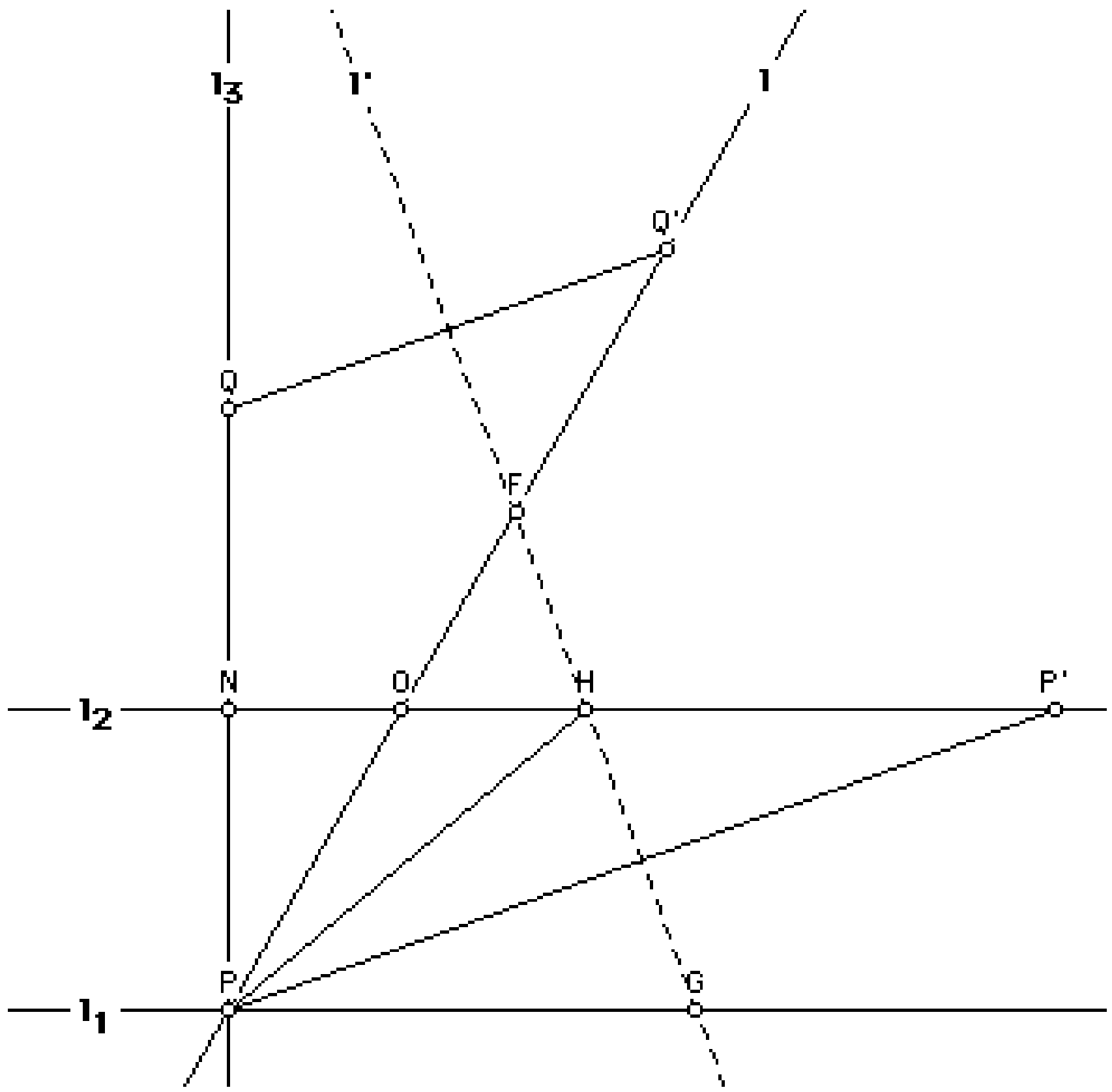

Figure 3: Abe's trisection.

The construction of a regular $n$-gon requires the construction of the real subfield of $\mathbb{Q}(\zeta_n)$, where $\zeta_n$ is a primitive $n$th root of unity. The real subfield of $\mathbb{Q}(\zeta_n)$ is $\mathbb{Q}(\cos(2\pi/n))$, a totally real field. It follows that $\zeta_n$ is constructible by bisections and trisections if and only if it is an element of $\mathbb{V} \oplus i\mathbb{V}$. Since $\mathbb{Q}(\zeta_n)$ has degree $\Phi(n)$ ($\Phi$ is Euler's totient function) over the field $\mathbb{Q}$ we have the following consequence of Theorem 3.1:

**Corollary 3.2.** *A regular $n$-gon can be constructed by origami using axioms (P), (L), (B), and (T) if and only if any prime factor of $\Phi(n)$ is either 2 or 3.*

The fundamental connection between trisection constructions and axiom (T) is given by:

**Theorem 3.3.** *Over a field that is closed under application of the Pythagorean axioms (P), (L), and (B) then any folded line $l'$ of axiom (T) can be constructed if and only if the angle between constructed lines can be trisected.*

*Proof.* Given axiom (T), we show how to trisect the angle between given (fold) lines $l$ and $l_1$ that meet at $P$. Since we can do all Pythagorean constructions, we can construct the perpendicular $l_3$ to $l_1$ at $P$. We next construct some point $Q$ on $l_3$ and then construct the perpendicular bisector $l_2$ of segment $PQ$. The line $l_2$ meets $l_3$ at $N$. Using axiom (T), we construct the line $l'$ such that the reflections $Q'$ and $P'$ of $Q$ and $P$ lie on $l$ and $l_2$, respectively. As shown in Figure 3, this line meets $l$ at $F$, $l_1$ at $G$, and $l_2$ at $H$. Then $\angle OPH$ is one-third of $\angle OPG$. To see this, let angle $\angle QP'N = \alpha$. Then $\angle NP'P = \alpha$, because $P'N$ is the altitude of the isosceles triangle $\triangle QP'P$. Since $NP'$ is parallel to $PG$, $\angle P'PG = \alpha$ as well. Using the isosceles triangle $\triangle PHP'$, we see also that $\angle HPP' = \alpha$. Now with the aid of the isosceles triangle $\triangle FP'P$, we observe that $2\alpha = \angle FP'P = \angle FPP'$, so $\angle OPG = 3\alpha$ and $\angle OPH = \alpha$

Conversely, suppose we can perform a trisection. Given (constructed) points $P$ and $Q$ and a (fold) line $l$ through $P$, as in Figure 3 ( $l$ is $OP$), we fold the perpendicular bisector of segment $PQ$ and the line $l_1$ that is perpendicular to $PQ$ through $P$. Now we trisect the acute angle formed by $l_1$ and $l$. The first of these trisection lines (one-third angle) intersects the line $l_2$ at a point $P'$, the other trisector (two-thirds angle) meets $l_2$ at a point $H$. The line through $Q$ and parallel to $PP'$ meets $l$ at a point $Q'$. If $F$ is the point where $l$ and the line $PQ'$ intersect, then (using congruent isosceles triangles) one can easily show that $FH$ is the fold line of axiom (T) as desired. □

**Folding the regular heptagon.** A method for construction of the regular heptagon using origami folding operations involving axiom (T) and Euclidean constructions is given in [2]. We indicate here how this can also be done using bisections and trisections. The first main step in the construction of the regular heptagon is a trisection that comes from the the extension $\mathbb{Q}(\cos(2\pi/7))$. The cubic equation associated with the heptagon is obtained by substituting $t = x + x^{-1}$ in

$$x^6 + x^5 + x^4 + x^3 + x^2 + x + 1 = 0$$

yielding

$$t^3 + t^2 - 2t - 1 = 0.$$

We first reduce this to

$$z^3 - \frac{7}{3}z - \frac{7}{27} = 0$$

by a change of variable. After a further substitution $u = z - 1/3$ and a scale change $w = \lambda u$ we arrive at the equation

$$4w^3 - 3w = 1/\sqrt{28}.$$

It follows that $w = \cos(\theta)$, where $\cos(3\theta) = 1/\sqrt{28}$ gives a solution.

Now to complete this heptagon with Pythagorean constructions requires a bisection or sum of squares expression for $\sin^2(2\pi/7)$. The Artin-Schreier theory predicts that a formula does exist. Here's one such:

$$4\sin^2(2\pi/7) = 3\cos^2(\pi/7) + 4\cos^4(3\pi/7).$$

# 4 TOTALLY REAL FIELDS OF DEGREE $2^a3^b$.

For a quadratic equation $x^2 + bx + c = 0$, the square of the difference of the roots (i.e., the discriminant) $\Delta_2 = (r_1 - r_2)^2 = (r_1 + r_2)^2 - 4r_1r_2 = b^2 - 4c$ is prominant in the algebraic formula that gives the roots. In solving a cubic $x^3 + Ax^2 + Bx + C = 0$ one can first simplify the equation by eliminating the coefficient of $x^2$. This operation, which is the analogue of completing the square in the quadratic case, just involves rewriting the equation in terms of $y = x - A/3$. This operation does not change the discriminant $\Delta_3 = (r_1 - r_2)^2(r_1 - r_3)^2(r_2 - r_3)^2$ of the cubic. A cubic has distinct roots if and only if $\Delta_3 \neq 0$. After completing the cubic, we obtain the reduced cubic $x^3 + px + q = 0$, whose discriminant is $\Delta_3 = -(27q^2 + 4p^3)$.

**Lemma 4.1.** *The cubic polynomial $f(x) = x^3 + px + q$ with coefficients in the real field $K$ has three distinct real roots if and only if $\Delta_3 > 0$. If $f(x)$ has three distinct real roots, then $p < 0$ and the cubic equation $x^3 + px + q = 0$ can be solved by a trisection over $K(\sqrt{3}, \sqrt{-p})$.*

*Proof.* It is easy to see that $f(x) = x^3 + px + q$ has three real roots if and only if this function has two critical points with critical values of opposite signs. The critical points occur at $\pm\sqrt{-p/3}$ and the product of the critical values is easily seen to be $-\Delta_3/27$. Thus there are three real roots if and only if $\Delta_3 \geq 0$.

Certainly if $p \geq 0$, then $\Delta_3 \leq 0$. Consequently, if there are three distinct roots (which requires that then $\Delta_3 > 0$), then $p$ must be negative. The classical trigonometric solution is obtained by making the substitution $x = y\sqrt{-4p/3}$. After multiplication by a suitable constant, the equation simplifies to $4y^3 - 3y = u$, where $|u| \leq 1$; in fact, $u^2 = -27q^2/4p^3$. Since $4\cos^3(\theta) - 3\cos(\theta) = \cos(3\theta)$ we can solve the cubic using trigonometry. The three roots correspond to $y = \cos(\theta)$ where $\theta$ is one of the three angles with $\cos(3\theta) = u$. If we take $\theta = \frac{1}{3}\mathrm{Arccos}(u)$, the other two angles are $\theta \pm \frac{\pi}{3}$. □

Thus, a cubic with three real roots can be solved by a trisection, after using a square root to obtain $u$. This method was certainly known in the sixteenth century to Viéte, who studied the formulas for expressing the sines and cosines of multiples of $\theta$ as polynomials in $\theta$.

**Totally real cubics.** Recall that the splitting field $L$ of an irreducible polynomial $f(x)$ of degree $n$ over a field $K$ is the smallest extension field of $K$ containing the roots of $f(x)$, i.e., $L = K(r_1, r_2, \ldots, r_n)$ is generated over $K$ by the roots $r_i$, $i = 1, 2, \ldots, n$ of $f(x) = 0$.

**Theorem 4.2.** *An irreducible cubic polynomial $f(x)$ over a totally real field $K$ has a totally real splitting field $L$ if and only if the roots of $f(x)$ can be construced using only field operations, extraction of square roots of totally positive elements in $K$, and angle trisection. In this case, the three roots of the cubic belong to a totally real extension of $K$ of degree which is either 3 or 6.*

*Proof.* By a change of variable, we may suppose that $f(x) = x^3 + px + q$ is a reduced irreducible cubic polynomial with coefficients belonging to the totally real field $K$. If the polynomial has three real roots, then necessarily $-p > 0$ and the splitting field $L$ over $K$ is a real field.

Under any embedding $\sigma : L \to \mathbb{C}$ the field $K$ is mapped into $\mathbb{R}$ since it is totally real. If $L$ is totally real then via this embedding the polynomial $f(x)$ is transformed to another polynomial with three real roots; hence the image of $-p$ under any isomorphism of $K$ into $\mathbb{R}$ is also positive otherwise the splitting field of the transformed polynomial will not be a real field. Thus $-p$ is a totally positve element of $K$. Thus as we have seen in Lemma 4.1, the three roots of $f(x)$ are obtained by a single trisection operation after adjoining to $K$ a square root of some totally positive element.

Conversely, suppose that one root of a cubic polynomial $f(x)$ with coefficients in $K$ can be constructed using the trisection of an angle and the extraction of square roots of totally positive elements of $K$, then the three (algebraic conjugate) numbers obtained from the trisection are the three real roots of $f(x)$. Hence the field generated over $K(\sqrt{-p}, \sqrt{3})$ by these three roots is totally real since it is real field and a splitting field, so invariant under embeddings into $\mathbb{C}$.

Moreover, since $f(x)$ is a cubic the splitting field $L$ of $f(x)$ is an extension over $K$ of degree which is either 3 or 6, depending on whether or not the discriminant is a square in $K$. $\square$

The next theorem expresses the consequences of our discussion in terms of towers of totally real fields of degree $2^a3^b$. Its proof, which is immediate from the foregoing results is similar to the tower theory for Euclidean or Mira constructions [5],[9].

**Theorem 4.3.** *Any real number $\alpha$ obtained from the rationals using field operations and a sequence of bisections or trisections is totally real and has a splitting field which is an extension of $\mathbb{Q}$ of degree $2^a3^b$ for some nonnegative integers $a$ and $b$.*

*Conversely, any totally real algebraic number $\alpha$ whose splitting field over the rationals has degree $2^a3^b$ is contained in the top field of a tower $\mathbb{Q} = F_0 \subset F_1 \subset \cdots \subset F_n$, where each extension either is Pythagorean of degree 2 or is of degree 3 and obtained by adjoining to the previous field a root of a trisection polynomial $4x^3 - 3x - u$, $|u| \leq 1$.*

The totally real subfield of the Viéten numbers $\mathbb{V}$ is the field denoted $\mathbb{BT}$.

**Corollary 4.4.** *The numbers constructible using axioms (P), (L), (B), and (T) are the elements of $\mathbb{BT} \oplus i\mathbb{BT}$.*

# 5 TRISECTIONS AND EUCLIDEAN CONSTRUCTIONS.

Adding trisections to the field of Euclidean constructions gives rise to a family of constructions studied by Gleason [6]. We may carry out all these constructions by application of the origami axioms using (P), (L), and (T) and the additional Euclidean axiom:

(E) For a given fold line $l$ and constructed points $P$ and $Q$ then (when possible) any line $l'$ through $Q$ with the property that the reflection of $P$ in $l'$ lies on $l$ can be constructed (folded).

As discussed in [1] and [2], the Euclidean construction axiom (E) allows for the extraction of all real square roots and thus includes even more constructions than the Pythagorean axiom (B). In fact the complex numbers which can be constructed by application of axioms (P), (L), and (E) is the field of Euclidean numbers [1].

The field obtained by constructions using the axioms (P), (L), (E), and (T) is denoted $\mathbb{ET}$. Of course, $\mathbb{ET}$ contains $\mathbb{BT}$, the field of totally real origami numbers. It is clear that any element constructed lies in a Galois (or normal) field extension generated over $\mathbb{Q}$ by quadratic extensions and totally real cubic extensions (by Theorem 4.2). It follows from the next theorem that $\sqrt[3]{2} \notin \mathbb{ET}$.

If an irreducible real polynomial $f(x)$ has only real roots we shall call it a totally real polynomial.

**Theorem 5.1.** *Any root $\beta$ of an irreducible cubic polynomial $f(x)$ with coefficients in $\mathbb{ET}$ that is not totally real can not be constructed using Euclidean tools and trisections.*

*Proof.* By way of contradiction we first treat the case of a real root. Suppose that the real number $\beta$ satisfies an irreducible real cubic equation that is not totally real but that we can construct $\beta$ using Euclidean tools and trisections. Equivalently assume that we can construct an extension $K$ of $\mathbb{Q}$ containing $\beta$ as the top field of a tower of quadratic or totally real cubic extensions. Let $F = E(\alpha)$ be the last cubic extension (generated by $\alpha$, the root of a real normal cubic) in this tower, so that the degree of $K$ over $F$ is a power of 2. Since $\beta$ belongs to $K$ and its degree over $F$ is odd, it follows immediately $\beta$ is a member of $F$. Thus we may assume that $F = K$. Let $R$ be the real subfield of $E$ obtained as the elements fixed under complex conjugation ($E$ is an extension of degree at most 2 over $R$). The extension $R(\alpha)$ is a splitting field for a totally real cubic, and $R(\alpha)$ is the real subfield of $E(\alpha)$. Because $\beta$ lies in $R(\alpha)$, a splitting field, its conjugates also lie in that field. On the other hand, by assumption not all the conjugates of $\beta$ are real. This contradiction shows that $\beta$ can not be constructed using Euclidean tools and trisections.

Suppose next that two of the roots, call them $\alpha$ and $\beta$ are nonreal. We may assume $f(x)$ is in reduced form: $f(x) = x^3 + px + q$. If both $\alpha$ and $\beta$ are in $\mathbb{ET}$, then the remaining (real) root is in $\mathbb{ET}$ as well, since the sum of the three roots is zero. However, that is impossible by the previous case, hence we will have at most one of the nonreal roots is in $\mathbb{ET}$. However, if one of these complex roots, say $\alpha$, is in $\mathbb{ET}$ then so is the other $\beta$, since $\beta$ satisfies a quadratic equation, $(\alpha + \beta)\alpha\beta = q$, with coefficients expressed in terms of $\alpha$ and $q$. This final contradiction finishes the proof of the proposition. □

Figure 4: Alhazen Pencil

Consider the real number $u = \sqrt{2 - \sqrt{2}}$. The field $\mathbb{Q}(u)$ is Euclidean constructible. Since $|u| \leq 1$, we can construct a solution $\alpha$ to $4x^3 - 3x - u = 0$ by a trisection. Let $E = \mathbb{Q}(u, \alpha)$. Because the conjugate of $u$, $u' = \sqrt{2 + \sqrt{2}}$, is greater than 1, the cubic polynomial $f(x) = 4x^3 - 3x - u'$ is not totally real. Let $\alpha'$ be the real root of $f(x)$. According to Theorem 5.1, the conjugate field $E' = \mathbb{Q}(u', \alpha')$ cannot be constructed with Euclidean tools and trisections. We record this observation as:

**Corollary 5.2.** *The Galois (normal) closure of a subfield of* $\mathbb{ET}$ *need not be a subfield of* $\mathbb{ET}$*.*

**Solving Alhazen's problem using a trisection.** Alhazen's problem of planar spherical optics can be formulated as: *Given points $a$ and $b$ in the plane exterior to the unit circle $\Gamma$ centered at the origin, locate $z$ on $\Gamma$ so that the angle $\angle azb$ is bisected by the diameter through $z$.*

As we now show the Alhazen problem admits a solution using Euclidean tools and trisections. Alhazan solved this problem in the tenth century A.D. by using a neusis construction, [11]. This appears to be quite similar to the trisection method originally due to Archimedes. However, the explanation of why Alhazen's problem can be solved with the aid of a trisection is a bit deeper. We begin with a review of the method presented in [2]. Regard the given points $a$ and $b$ as complex numbers. If $\arg(z)$ denotes the (principal) *argument* of a complex number, $\arg(\frac{a-z}{b-z})$ is the measure of the angle $\angle azb$. This gives an equation

$$\mathrm{Im}((ab)\overline{z}^2) = \mathrm{Im}((a+b)\overline{z}) \qquad (1)$$

that must be satisfied by any solutions $z$. It is easy to see that equation (1) is satisfied by $z = 1/\overline{a}$ and $z = 1/\overline{b}$.

Performing a preliminary rotation, we may assume that the positive $x$-axis bisects the angle $\angle aOb$. This ensures that $b = k\overline{a}$ for some $k > 0$ so $\mathrm{Im}(ab) = 0$. Let $q = \mathrm{Re}(ab)$, $r = \mathrm{Re}(a+b)$, $s = \mathrm{Im}(a+b)$. We can write equation (1) in terms of real coordinates as $2qxy + sx - ry = 0$, which is the equation of a hyperbola $H$. The midpoint of the segment $[1/\overline{a}, 1/\overline{b}]$ is the center of $H$. Consequently, the two solutions $1/\overline{a}$ and $1/\overline{b}$ are on different branches of $H$, and both are inside the unit circle. Hence $H$ meets $\Gamma$ in four distinct points. These four points of intersection can be solved using Euclidean methods and a trisector. A solution $z = x + iy$ to Alhazen's problem is a complex number such that $x$ and $y$ are simultaneous solutions to $A_1 = 0$ and $A_2 = 0$ where

$$A_1 = \ x^2 + y^2 - 1, \ A_2 = \ 2qxy + sx - ry.$$

These simultaneous solutions to $A_1 = 0$ and $A_2 = 0$ also lie on any conic in the Alhazen pencil $A_2 + \lambda A_1$. The degenerate conics in this pencil occur for the solutions $\lambda$ that are the roots of the cubic polynomial $f(\lambda) = \lambda^3 + \tau\lambda + qrs/2$, where $\tau = (s^2 + r^2 - 4q^2)/4$ (see [1]). Each of these degenerate conics consists of a pair of lines.

To show that $f(\lambda)$ has three real roots we notice first that there are four points at which $A_1 = 0$ and $A_2 = 0$ intersect. Joining the four points in pairs by lines gives six lines; pairing the six lines in three pairs so that each pair passes through all four points yields the three degenerate conics. Thus $f(\lambda)$ has three real roots.

Conversely, in terms of the algebra, a root of $f(\lambda)$ picks out a degenerate conic equation $A_2 + \lambda A_1 = 0$ in the pencil. We can factor that equation into lines by first factoring the highest degree (two) terms in the associated polynomial of $x$ and $y$. This involves solving quadratic equations. Once we have the six lines corresponding to the three roots of $f(\lambda)$, we can intersect them to get the four solutions to Alhazen's problem. There are three lines of the six that pass through each of the four solutions (see Figure 4).

We therefore have:

**Theorem 5.3.** *Alhazan's problem can be solved by Euclidean tools and trisections or equivalently, by using the origami constructions governed by (P), (L), (E), and (T).*

**ROGER C. ALPERIN** grew up in Miami, Florida, graduated from The University of Chicago and then received a Ph.D. from Rice University. His main areas of interest in mathematics involve the interplay between algebra and low-dimensional topology or geometry.

Department of Mathematics

San Jose State University
San Jose, CA 95192
alperin@math.sjsu.edu